\RequirePackage{fix-cm}
\documentclass{svjour3ao}                     

\smartqed  

\usepackage{graphicx, subfigure}
\usepackage{algorithm, algorithmic}
\usepackage{array}
\usepackage{amssymb,amsmath,graphics,pstricks,multirow}
\usepackage{mathtools}
\usepackage{tikz}
\usetikzlibrary{shapes.geometric}

\newcommand{\set}[1]{\left\lbrace #1\right\rbrace }
\newcommand{\RR}{\mathbb{R}}
\newcommand{\bnd}{\partial \Omega}

\newcolumntype{L}[1]{>{\raggedright\let\newline\\\arraybackslash\hspace{0pt}}m{#1}}
\newcolumntype{C}[1]{>{\centering\let\newline\\\arraybackslash\hspace{0pt}}m{#1}}
\newcolumntype{R}[1]{>{\raggedleft\let\newline\\\arraybackslash\hspace{0pt}}m{#1}}
\newcommand{\formExp}[1]{\cdot 10^{#1}}
\newcommand{\ssp}{\rule[0mm]{0mm}{9.5pt}}
\newcommand{\msp}{\rule[0mm]{0mm}{10.5pt}}

\usepackage{color}
\definecolor{marin}{rgb}   {0.,   0.3,   0.7}
\definecolor{rouge}{rgb}   {0.8,   0.,   0.}
\definecolor{sepia}{rgb}   {0.8,   0.5,   0.}
\usepackage[colorlinks,citecolor=marin,linkcolor=rouge,
            bookmarksopen,
            bookmarksnumbered
           ]{hyperref}

\tikzset{
	triangle/.style={
		draw,
		shape border rotate=0,
		regular polygon,
		regular polygon sides=3,
		fill=black,
		node distance=2cm,
		minimum height=0.05cm,
		minimum width=0.05cm,
		text width=0.01cm,
		thin,
		inner sep=1pt
	}
}

\date{Version of November 6, 2017}

\journalname{}
\begin{document}

\title{Efficient boundary corrected Strang splitting}
\author{Lukas Einkemmer \and Martina Moccaldi \and Alexander Ostermann}
\institute{L.~Einkemmer\at Department of Mathematics, University of Innsbruck, 6020 Innsbruck, Austria\\ \email{lukas.einkemmer@uibk.ac.at} \and
M.~Moccaldi\at  Department of Mathematics, University of Salerno, Fisciano (SA), Italy\\ \email{mmoccaldi@unisa.it} \and
A.~Ostermann \at Department of Mathematics, University of Innsbruck, 6020 Innsbruck, Austria\\ \email{alexander.ostermann@uibk.ac.at}}
\maketitle

\begin{abstract}
Strang splitting is a well established tool for the numerical integration of evolution equations. It allows the application of tailored integrators for different parts of the vector field. However, it is also prone to order reduction in the case of non-trivial boundary conditions. This order reduction can be remedied by correcting the boundary values of the intermediate splitting step. In this paper, three different approaches for constructing such a correction in the case of inhomogeneous Dirichlet, Neumann, and mixed boundary conditions are presented. Numerical examples that illustrate the effectivity and benefits of these corrections are included.
\end{abstract}



\section{Introduction} \label{sec:intro}
In recent years, splitting schemes have been widely used to integrate evolution equations (see, for instance, \cite{Dawson1992,Estep2008,Hundsdorfer2003}). The main advantage of the splitting approach, compared to standard integrators, is the fact that it allows the separate treatment of different terms of the vector field. This usually reduces the computational effort required to solve the problem. Moreover, this strategy facilitates in many cases a parallel implementation and may lead to numerical methods having better geometric properties. For all these reasons, splitting methods are highly suitable for the numerical treatment of models describing complex phenomena. For instance, in turbulent combustion simulations with finite-rate kinetics, hundreds of differential equations (one for each chemical species) have to be integrated in addition to the Navier--Stokes equations for momentum and energy. Moreover, such problems are typically stiff, so they require implicit solvers. In such examples, splitting separately treats the  chemical kinetics and the convection-diffusion terms. This allows one to use appropriate numerical methods for each part and can thus considerably reduce the computational cost \cite{Stone2016}.

This paper is focused on diffusion-reaction systems. Such problems are efficiently integrated by splitting schemes \cite{Descombes2001,Gerish2002}. In particular, we will consider equations where the diffusion is modelled by a linear elliptic differential operator and the reaction by a nonlinear function. The splitting approach takes into account the different nature of the components and replaces the nonlinear system of differential equations with a linear system and a set of ordinary differential equations, making the implementation easier and reducing the computational burden. Moreover, splitting methods preserve positivity and invariant sets provided that the integrator of each flow has the corresponding property \cite{Hansen2012}.

When periodic or homogeneous Dirichlet boundary conditions are employed, the well-known Strang splitting scheme is second-order accurate, in general. However, it has been proved  in the literature that the order of Strang splitting is reduced in case of non trivial boundary conditions, e.g. inhomogeneous Dirichlet, Neumann, Robin boundary conditions or mixed ones (see, for instance, \cite{EOPart1,EOPart2,Hundsdorfer2003}). One possibility to overcome this order reduction is introducing a smooth correction function in splitting schemes such that the new reaction flow is compatible with the prescribed boundary conditions \cite{EOPart1,EOPart2}. For time-invariant Dirichlet boundary conditions, this correction has to be computed only once at the beginning of the simulation. However, for time-dependent Dirichlet, Neumann, or Robin boundary conditions, the correction is time-dependent. In these cases, the correction function has to be computed at each time step, which could noticeably increase the computational cost. Hence, an open problem is finding an efficient procedure to construct this function in order to preserve the aforementioned computational advantages of the splitting scheme.

In this work, we present and compare different strategies to construct the correction function. First of all, we propose to select the correction function as the solution of an elliptic problem endowed with appropriate boundary conditions \cite{EOPart1,EOPart2}. A similar elliptic problem has to be solved by the splitting method in each time step. Therefore, the computational overhead of this additional solve is moderate. Nevertheless, we investigate here still other approaches which, depending on the situation, are less expensive. E.g., an effective procedure consists in computing the correction function with the help of the actual numerical solution. For Dirichlet boundary conditions, the correction is computed at each time step as $f(u_n)$, where $f$ is the vector field of the reaction and $u_n$ is the numerical solution at the beginning of the current time step. This is a cheap option for problems where the application of $f$ is cheap. However, there are situations in which additional evaluations of the vector field $f$ highly increase the computational effort, as it is the case in turbulent combustion simulations \cite{Stone2016}. Therefore, we present two further techniques which are more convenient for problems having expensive reactions. For Dirichlet boundary conditions, we employ a widely used low-pass filter for noise reduction in linear image processing \cite{Acharya2005,Smith1999}, the moving average, which we extend with multiple grid levels to further increase smoothness. In case of Neumann boundary conditions, we use a crude approximation of the elliptic problem on a hierarchy of grids, which results in a multigrid-like algorithm.

The outline of the paper is as follows. In Section \ref{sec:problem} we present the splitting approach, explain the related problem of order reduction and the use of a correction function to avoid it. Sections \ref{sec:Dirichlet} and \ref{sec:Neumann} are devoted to the description of efficient algorithms to construct the correction function for Dirichlet and Neumann boundary conditions, respectively. Moreover, we report some numerical tests that prove the effectiveness of the proposed approaches. Finally, we discuss some conclusions in Section \ref{sec:Conclusions}.


\section{Model problem} \label{sec:problem}
In this paper, we deal with the numerical integration of diffusion-reaction systems modelled by the following initial-boundary value problem
\begin{subequations} \label{system}
\begin{alignat}{3}
\partial_t u &= D u + f(u),  \label{equation}\\
B u|_{\partial \Omega} &= b, \label{BCs}\\
u(0)&=u_0, \label{ICs}
\end{alignat}
\end{subequations}
where $u: [0,T] \times \Omega \to \RR $. Here, $\Omega \subset \RR^d$ is a bounded domain with sufficiently smooth boundary $\partial \Omega$, the reaction term $f$ is a real smooth function, typically nonlinear. Moreover, $D$ is an elliptic differential operator (e.g. the Laplacian), the boundary data $b:[0,T] \times \partial \Omega$ may be time-dependent and
\begin{equation}\label{Boperator}
	B=  \beta(x)\partial_n + \alpha(x)
\end{equation}
is a first-order differential operator with sufficiently smooth coefficients. The functions $b$ and $u_0$ are also assumed to be sufficiently smooth. According to the values of the coefficients in $B$, different boundary conditions can be modelled. For example, if $\beta(x)=0 $ and $ \alpha(x)\ge c >0$, then \eqref{BCs} corresponds to Dirichlet boundary conditions, whereas $\alpha(x)=0$ and $\beta(x)\ge c >0$ give Neumann boundary conditions.

For the numerical solution of \eqref{system} we split the system into the linear diffusion equation
\begin{subequations}\label{Prob12}
\begin{equation}\label{Prob1}
\partial_t v = D v, \quad B v |_{\partial \Omega} = b
\end{equation}
and the nonlinear reaction equation
\begin{equation} \label{Prob2}
\partial_t w  = f(w).
\end{equation}
\end{subequations}
Let $u_n$ be the numerical approximation to the exact solution $u$ of \eqref{system} at time $t=t_n$. To step from $t_n$ to $t_{n+1}=t_n+ \tau$, where $\tau$ is the step size, we use the well-known Strang splitting approach
\begin{equation} \label{Strang}
	u_{n+1}=\mathcal{S}_{\tau} u_n = \varphi_{\tau/2}(\psi_{\tau}(\varphi_{\tau/2}(u_n))),
\end{equation}
where $\varphi_{\tau} $ and $\psi_{\tau}$ are the exact flows of the first \eqref{Prob1} and second \eqref{Prob2} sub-problem, respectively. In other words, Strang splitting for the integration of problem \eqref{system} consists in computing the solution $v(\frac{\tau}{2})$ of \eqref{Prob1} with initial value $v(0)=u_n$, then integrating problem \eqref{Prob2} with initial value $w(0)=v(\frac{\tau}{2})$ to obtain $w(\tau)  $ and finally computing the solution of \eqref{Prob1} with initial value $v(\frac{\tau}{2}) = w(\tau)  $ to get $v (\tau) = u_{n+1} $.

For diffusion-reaction equations, a convergence proof for Strang splitting was presented in \cite{Hansen2012}. There it was shown that second-order convergence requires a certain consistency between the boundary conditions of the diffusion problem and the action of the reaction flow. In general, however, this assumption is not satisfied for inhomogeneous boundary conditions and order reduction takes place. To avoid this problem, the authors of \cite{EOPart1,EOPart2} propose to employ a correction function $q$ which is sufficiently smooth and fulfills the boundary conditions of the nonlinearity, i.e.
\begin{equation}
	 B q |_{\partial \Omega} =  B f(u) |_{\partial \Omega}.
\end{equation}
With the help of $q$ they modify system \eqref{Prob12} as follows:
\begin{subequations}\label{ModSplitSystem}
	\begin{alignat}{2}
	&\partial_t v = D v+q, \quad B v |_{\partial \Omega} = b, \label{ModProb1}\\
	& \partial_t w  = f(w)-q. \label{ModProb2}
	\end{alignat}
\end{subequations}
An application of Strang splitting to this system is called \emph{modified} Strang splitting in the following presentation and has been shown to be extremely competitive compared to other boundary correction techniques \cite{einkemmer2017dcssb}. We remark that this approach, in general, requires that the correction function is computed once every time step.

In \cite{EOPart1,EOPart2} the correction function has been computed solving the corresponding elliptic problem. However, this strategy increases the computational cost. Indeed, in case of time-invariant Dirichlet boundary conditions, the correction function can be precomputed since it does not rely on the current time step and the solution $u$. For time-dependent Dirichlet and Neumann boundary conditions, however, the correction function does depend on time and in the latter case on the numerical solution $u$. Therefore, it has to be updated at each time step, which increases the computational burden. Integrating the elliptic problem is not the only way to compute the correction function, as will be shown in this paper. Also note that the choice of the correction function has some influence on the accuracy of the method \cite{EOPart2}. Therefore, finding a good and efficient strategy to construct this function is crucial to overcome the order reduction and to obtain an efficient splitting scheme. In this treatise, we present various techniques to deal with this problem.


\section{The Dirichlet case} \label{sec:Dirichlet}

Let us consider system \eqref{system} equipped with inhomogeneous Dirichlet boundary conditions
\begin{equation}\label{DBCs}
	u |_{\partial \Omega} = b.
\end{equation}
Motivated by the investigation carried out in \cite{EOPart1}, we aim to construct a correction function $q$ such that
\begin{equation} \label{DCorrFunc}
	q |_{\partial \Omega} = f(b).
\end{equation}
For the modified Strang splitting, it is sufficient to know $q$ at times $t_n= n\tau$. We denote these approximations by $q_n$ and set $b_n=b(t_n)$.
In the following, we present three different strategies to achieve this goal.

A first attempt consists in constructing the correction function by solving the elliptic problem
\begin{equation}\label{SolveProblem}
\begin{split}
D q_n &= 0, \\
q_n|_{\partial \Omega} &= f(b_n).
\end{split}
\end{equation}
The integration of this system can be carried out, for instance, through a direct method (as in numerical tests presented in this paper) or through Krylov subspace methods. However, this strategy could be more expensive than the other procedures presented below. We also recall that $q$ only has to satisfy \eqref{DCorrFunc}. Therefore, the differential operator $D$ in \eqref{SolveProblem} can be replaced by any other elliptic differential operator (e.g.~the Laplacian $\Delta$) for which the problem can be solved efficiently.

The second approach for finding a correction $q_n$ makes direct use of the numerical solution. Since the numerical solution $u_n$ satisfies the boundary condition at $t=t_n$, we can choose the correction function as
\begin{equation}\label{Ddirect}
	q_n=f(u_n).
\end{equation}
This procedure is easy to implement and leads to a smooth correction function under the assumption that the reaction is smooth.  Moreover, it does not depend on the structure of the domain. However, evaluating $f$ at each time step could be computationally expensive, as is the case in the situation described in \cite{Stone2016}. Therefore, it can be worthwhile to consider an alternative approach.

The basic idea of this third strategy is to employ a very efficient spatial filtering technique used for noise reduction (or smoothing) in linear image processing (see, for instance, \cite{Acharya2005,Smith1999}): the moving average filtering or box blur. We recall that a filter is a device or process which removes some unwanted components or features from given data, e.g.~an input image. Suppose that the input image is composed of $R \times S$ pixels. In the case of moving average the output image is obtained by substituting each pixel in the input image with a weighted average of its neighbouring pixels, as follows:
\begin{equation}\label{outputImage}
g_{\mathrm{out}}(r,s) = \sum_{(i,j)\in (r,s) +\mathcal{W}} w_{i-r,j-s}\; g_{\mathrm{in}}(i,j), \qquad r=0,\dots,R-1, \ s=0,\dots, S-1,
\end{equation}
where $g_{\mathrm{in}}(r,s)$ is a pixel in the input image, $g_{\mathrm{out}}(r,s)$ the corresponding one in the output image, $\mathcal{W}$ is a neighbourhood of the considered pixel and $w_{i-r,j-s}$ are the filter weights.

For the construction of the correction function $q_n$, however, applying the moving average on the finest grid does not yield sufficient smoothness. For this reason, we propose to use the moving average filtering in a multigrid-like fashion.

Multigrid schemes are commonly employed for the iterative solution of elliptic problems. They scale linearly with the number of unknowns and thus significantly accelerate the convergence of the basic iterative method \cite{Briggs1987,Hackbusch1985}. The main observation is that low frequencies on a fine grid become high frequencies on a coarse grid. The multigrid strategy consists in relaxing the linear system (derived from the discretization of the considered PDE) on a hierarchy of grids. As a consequence, such methods require mechanisms to transfer information from a coarse grid to a fine grid (\emph{prolongation}) and, conversely, from a fine grid to a coarse one (\emph{restriction}) (see Section \ref{sec:Neumann} for further details). In particular, we use a sort of one half of a V-cycle, from bottom to top, i.e. from the coarse grid to the fine grid. In this case, the prolongation is carried out through the moving average \eqref{outputImage}.

\begin{figure}
	\centering
	\begin{tabular}{C{5.9cm}C{5.9cm}}
	\subfigure[]	
	{\begin{tikzpicture}
	\draw (0,0) rectangle (4,4);
	\draw [dashed] (0,2) -- (4,2);
	\draw [dashed] (2,0) -- (2,4);
	\node[right] at (4.2,4) {$\Omega$};
	\draw [<->] (4.3,0) -- (4.3,2);
	\node[right] at (4.31,1) {$2H$};
	\draw [fill=black] (0,0) circle [radius=0.07];
	\draw [fill=black] (2,0) circle [radius=0.07];
	\draw [fill=black] (4,0) circle [radius=0.07];
	\draw [fill=black] (0,2) circle [radius=0.07];
	\draw [fill=black] (2,2) circle [radius=0.07];
	\draw [fill=black] (4,2) circle [radius=0.07];
	\draw [fill=black] (0,4) circle [radius=0.07];
	\draw [fill=black] (2,4) circle [radius=0.07];
	\draw [fill=black] (4,4) circle [radius=0.07];

	\draw              (3,1) circle [radius=0.07];
	\node[above left] at (3,1) {$C_1$};
	\draw              (3,3) circle [radius=0.07];
	\node[above left] at (3,3) {$C_2$};
	\draw              (1,3) circle [radius=0.07];
	\node[above left] at (1,3) {$C_3$};
	\draw              (1,1) circle [radius=0.07];
	\node[above left] at (1,1) {$C_4$};
	\end{tikzpicture}\label{Fig:grid1}}
	&
	\subfigure[]	
	{\begin{tikzpicture}
	\draw (0,0) rectangle (4,4);
	\draw [dashed] (0,2) -- (4,2);
	\draw [dashed] (2,0) -- (2,4);
	\draw [dashed] (1,0) -- (1,4);
	\draw [dashed] (3,0) -- (3,4);
	\draw [dashed] (0,1) -- (4,1);
	\draw [dashed] (0,3) -- (4,3);
	\node[right] at (4.2,4) {$\Omega$};
	\draw [<->] (4.3,0) -- (4.3,1);
	\node[right] at (4.31,0.5) {$H$};
	\draw [fill=black] (3,3) circle [radius=0.07];
	\node[above left] at (3,3) {$C_2$};
	\draw [fill=black] (1,3) circle [radius=0.07];
	\node[above left] at (1,3) {$C_3$};
	\draw [fill=black] (1,1) circle [radius=0.07];
	\node[above left] at (1,1) {$C_4$};
	\draw [fill=black] (3,1) circle [radius=0.07];
	\node[above left] at (3,1) {$C_1$};
	\draw [fill=white](3,2) circle [radius=0.07];
	\node[above left] at (3,2) {$P_1$};
	\draw [fill=white](1,2) circle [radius=0.07];
	\draw [fill=white](2,3) circle [radius=0.07];
	\node[above left] at (2,3) {$P_2$};
	\node[above left] at (1,2) {$P_3$};
	\draw [fill=white](2,1) circle [radius=0.07];
	\node[above left] at (2,1) {$P_4$};
	\end{tikzpicture} \label{Fig:grid2}}\\
	\multicolumn{2}{c}{
	\subfigure[]{
	\begin{tikzpicture}
	\draw (0,0) rectangle (4,4);
	\draw [dashed] (0,2) -- (4,2);
	\draw [dashed] (2,0) -- (2,4);
	\draw [dashed] (1,0) -- (1,4);
	\draw [dashed] (3,0) -- (3,4);
	\draw [dashed] (0,1) -- (4,1);
	\draw [dashed] (0,3) -- (4,3);
	\node[right] at (4.2,4) {$\Omega$};
	\draw [<->] (4.3,0) -- (4.3,1);
	\node[right] at (4.31,0.5) {$H$};
	\draw [fill=black] (2,0) circle [radius=0.07];
	\draw [fill=black] (4,0) circle [radius=0.07];
	\draw [fill=black] (4,2) circle [radius=0.07];
	\draw [fill=black] (4,2) circle [radius=0.07];
	\draw [fill=black] (0,0) circle [radius=0.07];
	\draw [fill=black] (0,1) circle [radius=0.07];
	\draw [fill=black] (0,2) circle [radius=0.07];
	\draw [fill=black] (0,3) circle [radius=0.07];
	\draw [fill=black] (0,4) circle [radius=0.07];
	\draw [fill=black] (4,0) circle [radius=0.07];
	\draw [fill=black] (4,1) circle [radius=0.07];
	\draw [fill=black] (4,2) circle [radius=0.07];
	\draw [fill=black] (4,3) circle [radius=0.07];
	\draw [fill=black] (4,4) circle [radius=0.07];
	\draw [fill=black] (1,4) circle [radius=0.07];
	\draw [fill=black] (2,4) circle [radius=0.07];
	\draw [fill=black] (3,4) circle [radius=0.07];
	\draw [fill=black] (1,0) circle [radius=0.07];
	\draw [fill=black] (3,0) circle [radius=0.07];
	\node[triangle] at (2,2){ };
	\node[above left] at (3,3) {$C_2$};
	\node[triangle] at (3,3){ };
	\node[above left] at (1,3) {$C_3$};
	\node[triangle] at (1,3){ };
	\node[above left] at (1,1) {$C_4$};
	\node[triangle] at (1,1){ };
	\node[above left] at (3,1) {$C_1$};
	\node[triangle] at (3,1){ };
	\node[above left] at (3,2) {$P_1$};
	\node[triangle] at (3,2){ };
	\node[above left] at (2,3) {$P_2$};
	\node[triangle] at (2,3){ };
	\node[above left] at (1,2) {$P_3$};
	\node[triangle] at (1,2){ };
	\node[above left] at (2,1) {$P_4$};
	\node[triangle] at (2,1){ };
	\end{tikzpicture} \label{Fig:grid3}} }
	\end{tabular}
	\caption{Illustration of Algorithm \ref{algorithm:moving_average} for $k=2$. In (a) and (b) full points are used in the averages \eqref{q_movingAverage} to compute the values of $q^H$ in the empty points. In (c) the values of $q^H$ in the triangle-shaped points are updated by employing all direct neighbours, as described in step \ref{eq:step4} of Algorithm \ref{algorithm:moving_average}\black.\label{Fig:grid}}
\end{figure}
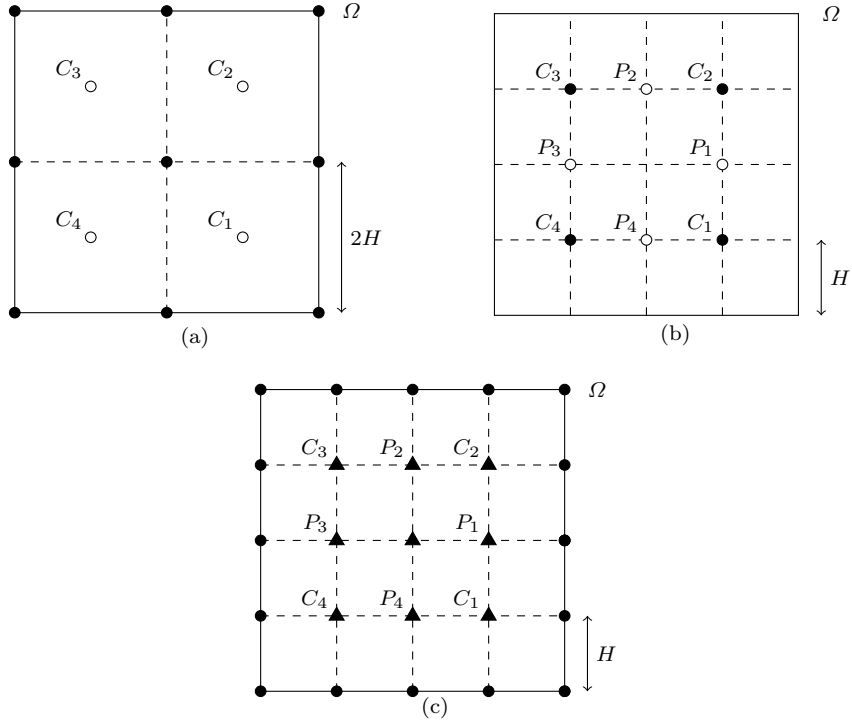

For the ease of presentation, we suppose that $d=2$ and $\Omega=(0,1) \times (0,1)$. However, this procedure can also be applied to more general domains. Following the described approach, we discretize $\overline\Omega$ in $(M+1)^2$ points with $M=2^\ell$ for simplicity. We thus consider the spatial mesh
\begin{equation*}
\Omega^\ell = \left\{(ih,jh) \mathrel{\big|} 0\le i,j\le 2^\ell \text{ and } h=2^{-\ell}\right\}
\end{equation*}
with mesh width $h$ representing the finest grid. We recall that the values of $q$ on the boundary $\bnd$ are known thanks to condition \eqref{DCorrFunc} and must not be changed. Therefore, we aim to compute the values of the correction function $q$ at each internal grid point.  For this purpose, we consider a hierarchy of grids $\Omega^k$ for $1\le k \le \ell$. Starting from the coarsest grid, the values of the correction function $q$ are computed at each internal grid point of the next finer grid $\Omega^k$ by using the moving average \eqref{outputImage}. The grid function $q^H$ with $H=2^{-k}$ is prescribed on the boundary $\partial\Omega$. For interpolating its values at the remaining (inner) points of $\Omega^k$, we employ the following relations, where $N=2^{k-1}$. We first set
\begin{subequations}\label{q_movingAverage}
\begin{equation}\label{eq:qmA1}
q^{H}(2r,2s) = q^{2H}(r,s)
\end{equation}
for $r,s=1,\dots,N-1$, then compute
\begin{equation}\label{eq:qmA2}
\begin{aligned}
q^{H}(2r+1,2s+1) &= \frac14 \Big( q^{2H}(r,s) + q^{2H}(r+1,s) \\
&\qquad\quad + q^{2H}(r,s+1) + q^{2H}(r+1,s+1)\Big)
\end{aligned}
\end{equation}
for $0\le r,s\le N-1$ and finally define
\begin{align}
q^{H}(2r+1,2s) &= \frac12 \Big( q^{H}(2r+1,2s-1) + q^{H}(2r+1,2s+1)  \Big),\label{eq:qmA3}\\
q^{H}(2r,2s+1) &= \frac12 \Big( q^{H}(2r-1,2s+1) + q^{H}(2r+1,2s+1) \Big)\label{eq:qmA4}
\end{align}
\end{subequations}
for the remaining inner points. At the end of this iteration, the following smoothing step is applied
\begin{equation}\label{SmoothingStep}
q^H(r,s) = \dfrac{1}{5} \sum_{(i,j)\in (r,s) +\mathcal{W}_5} q^{H}(i,j), \qquad r,s=1,\dots,2N-1,
\end{equation}
where $\mathcal{W}_5=\set{(0,0),\, (\pm 1, 0), \, (0, \pm 1)}$ is a neighbourhood of 0. The whole procedure for a single square subdomain is shown in Figure \ref{Fig:grid} and described in Algorithm~\ref{algorithm:moving_average}. Starting from the boundary values, this algorithm provides the values of the correction function $q(ih,jh)=q^h(i,j)$ on the finest grid $\Omega^h$.

\begin{figure}[htb]
 	\centering
	{\includegraphics[width=9cm]{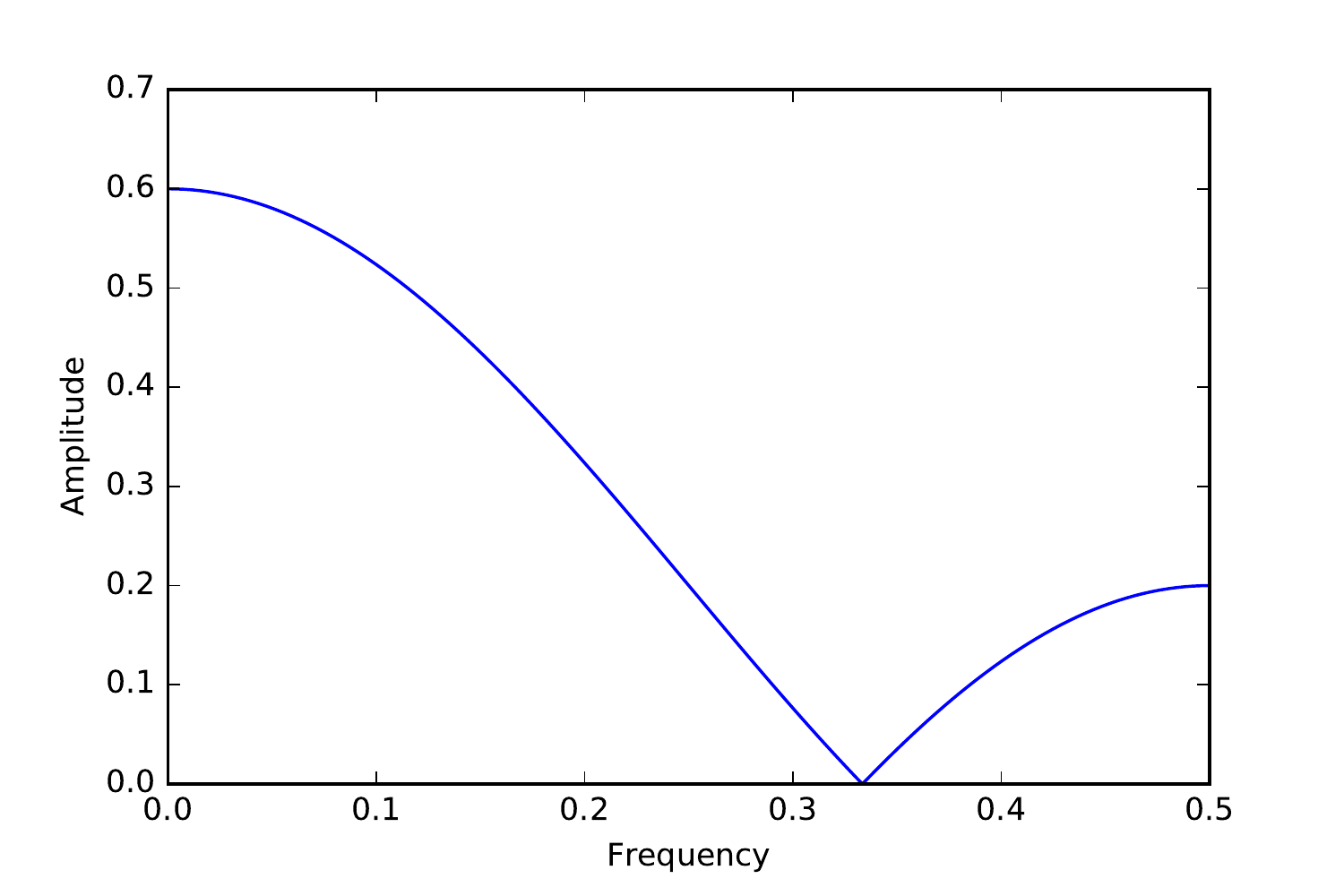}}
	\caption{Damping of frequencies due to the moving average approach \eqref{outputImage}. The plot has been generated by using the expression of the amplification factor \eqref{AmplFact}, fixing $ \frac{n}{M+1}=0.25$ and varying $\frac{m}{M+1}$ from $0$ to $0.5$.}\label{Fig:Fourier}
\end{figure}

Let us show that this approach leads to a smooth correction function $q$. In image processing, the moving average \eqref{outputImage} is a low-pass filter, i.e. a filter which returns an output image smoother than the original one and devoid of high spatial frequencies. This damping of high frequencies can be proved by converting the relations \eqref{q_movingAverage} and \eqref{SmoothingStep} from the space domain to the frequency domain through Fourier analysis (see, for instance \cite{Hackbusch1985,Smith1999} and reference therein). We illustrate this with formula \eqref{SmoothingStep} used in the smoothing step (i.e.~step \ref{eq:step4} of Algorithm \ref{algorithm:moving_average}). In the frequency domain it has the following representation:
\begin{equation}\label{MovingAverage_Fdomain}
	\widehat{Q}_{m,n}^{\text{new}} = \dfrac{1}{5}
	\left(1+\textrm{e}^{\textrm{i} 2 \pi\frac{m}{M+1}} +
	\textrm{e}^{-\textrm{i} 2 \pi\frac{m}{M+1}}+
	\textrm{e}^{\textrm{i} 2 \pi\frac{n}{M+1}}+
	\textrm{e}^{-\textrm{i} 2 \pi\frac{n}{M+1}}	
	\right)\widehat{Q}_{m,n},
\end{equation}
where $\widehat{Q}_{m,n}$  and $\widehat{Q}_{m,n}^{new}$ are the discrete Fourier transforms of $ q^H(r,s) $  before and after the smoothing step, respectively. Therefore, the factor of amplification of frequencies related to the moving average \eqref{SmoothingStep} at $(m,n)$ is the following:
\begin{equation} \label{AmplFact}
\rho(m,n)= \dfrac{1}{5} \left|1+2 \cos\left(  \dfrac{2 \, \pi \, m}{M+1}\right) +2 \cos\left(  \dfrac{2 \, \pi \, n}{M+1}\right) \right|.
\end{equation}
Figure \ref{Fig:Fourier} shows the low-pass feature of the moving average \eqref{SmoothingStep}. Indeed, some high frequencies are completely damped, whereas other high frequencies are attenuated by a factor of $1/10$ or $2/10$. This figure is obtained by fixing $\frac{n}{M+1}=0.25$ and varying $\frac{m}{M+1}$ from $0$ to $0.5$.

\begin{algorithm}[htb]
\caption{How to compute the values of the correction function $q$ in the square domain $\Omega =(0,1)^2$ through the moving average approach.}
\label{algorithm:moving_average}
\vspace{1mm}
\begin{enumerate}
\item Let $\Omega^0$ be the four corners of $\Omega$, define $q^1(i,j)$ for $0\le i,j\le 1$ by using the given boundary values and set $k=1$.
\vspace{1mm}
\item Set $H=2^{-k}$ and define $q^H$ at all even grid points by \eqref{eq:qmA1}.\label{eq:step2}
\item Cut the diagonals of the cells of $\Omega^{k-1}$ to get new centers $C_\kappa$ (as indicated in Figure \ref{Fig:grid1} for $k=2$).
\vspace{1mm}
\item Compute the values of $q^H$ at each center $C_\kappa$ by taking the average \eqref{eq:qmA2} among the neighbouring vertices of $\Omega^{k-1}$ (see Figure \ref{Fig:grid1}).
\vspace{1mm}
\item Use the average \eqref{eq:qmA3} (or the average \eqref{eq:qmA4}) to compute the values of $q^H$ at the midpoints $P_\kappa$ of each horizontal (or vertical) segment connecting two neighbouring centers $C_\mu$ and $C_\nu$ (see Figure \ref{Fig:grid2}).
\vspace{1mm}
\item Define the missing values of $q^H$ on $\partial\Omega\cap\Omega^k$ by using the given boundary values.
\vspace{1mm}
\item Recompute the values of $q^H$ at each interior point of the grid $\Omega^k$ by carrying out the \emph{smoothing step} \eqref{SmoothingStep}  (the triangle-shaped points in Figure \ref{Fig:grid3}).\label{eq:step4}
\vspace{1mm}
\item If $k<\ell$, set $k=k+1$ and go to step \ref{eq:step2}.
\vspace{1mm}
\item On the finest grid $\Omega^\ell$ where $h=2^{-\ell}$, set $q(ih,jh) = q^h(i,j)$ for $0\le i,j\le 2^\ell$.
\end{enumerate}
\end{algorithm}

\subsection{Numerical results} \label{subsec:DirichletNumExp}

In this section, we present some numerical tests for the following diffusion-reaction problem
\begin{subequations}\label{DR_Dproblem}
\begin{alignat}{3}
\partial_t u &= \Delta u + u^2, \label{DR_Dequation} \\
u|_{\partial \Omega} &= b, \label{DR_DBCs}\\
u(0)&=u_0 \label{DR_DICs}	
\end{alignat}
\end{subequations}
on $\Omega=(0,1)^2$. In all the presented experiments the Laplacian is discretized by the classic second-order centered finite differences and the resulting problem is integrated with the modified splitting approach described in Section \ref{sec:problem}. In particular, we employ Strang splitting and we solve each problem in \eqref{ModSplitSystem} by means of the variable step size variable order multistep scheme \textsf{VODE}, which is based on the backward differentiation formulas. As relative error tolerance, we prescribe $10^{-10}$. Moreover, we emphasize that all the numerical tests in this paper are conducted in the stiff regime (i.e. $\tau \gg h^2$, where $\tau$ is the time step size and $h$ is the spatial mesh width).

For a grid function $E$ on $\Omega^h$ with zero boundary values, we define the discrete $L^2$ norm (called $l^2$ in the tables below) \black and the $l^\infty$ norm in the usual way:
\begin{equation}\label{eq:norms}
\|E\|_2 = \sqrt{\frac1{(M-1)^2} \sum_{i,j=1}^{M-1} |E_{i,j}|^2},\qquad \|E\|_\infty = \max_{1\le i,j\le M-1} |E_{i,j}|.
\end{equation}

In the remainder of this section, we employ the expressions \emph{Strang} and \emph{modified Strang} to refer to the standard Strang splitting (i.e, \eqref{Strang} applied to \eqref{Prob12}) and the modified Strang splitting (i.e., \eqref{Strang} applied to \eqref{ModSplitSystem}), respectively. In order to distinguish the different strategies for computing the correction function, we use the following expressions. If the correction function is computed as the solution of the Dirichlet problem \eqref{SolveProblem}, we call it \emph{exact elliptic}. If it is computed as the moving average of neighbouring points in the way described in Algorithm~\ref{algorithm:moving_average}, we call it \emph{grid average}. Finally, the computation of $q$ by employing the values of the numerical solution at the beginning of the current time step \eqref{Ddirect} is labelled as $q_n=f(u_n)$.

\begin{example}
We now integrate problem \eqref{DR_Dproblem} with the following boundary and initial data:
\begin{subequations}\label{test1}
\begin{alignat}{2}
	b(t,x,y) &= 1+\dfrac{t}{1+x^2+y^2+t^2}, \qquad (x,y)\in\partial\Omega,\label{test1_DBCs}\\
	u(0,x,y) &= 1+ \sin^2 (\pi x) \, \sin^2 (\pi y),\qquad (x,y)\in \Omega. \label{test1_ICs}
\end{alignat}
\end{subequations}
Tables \ref{Tab:D_InfNorm} and \ref{Tab:D_2Norm} report the errors and the observed orders of convergence for all considered schemes. The errors are computed in the infinity (Table \ref{Tab:D_InfNorm}) and the $2$-norm (Table \ref{Tab:D_2Norm}) with respect to a reference solution obtained by solving the unsplit problem with \textsf{VODE} and prescribed tolerance equal to  $10^{-10}$. For standard Strang splitting, we observe the expected reduction \cite{EOPart1} to orders $1$ and $1.25$ in the infinity and the $2$-norm, respectively. The modified Strang splitting schemes, however, are always of order two. Therefore, they are all suitable to overcome the problem of order reduction. Moreover, they are all significantly more accurate than standard Strang splitting, even if the modified Strang splitting based on the direct construction $q_n=f(u_n)$ exhibits a slightly larger error in the infinity norm and a smaller error in $2$-norm than the other modified Strang splitting schemes. Finally, Figure \ref{Fig:DCorrFunc} shows that the  correction functions constructed by the proposed methods are different, but they are all smooth and they are all compatible with condition \eqref{DCorrFunc}, as required.

\begin{center}
	\begin{table}[h!t!b!]
		\caption{Errors and observed orders of standard and modified Strang splitting schemes for the integration of problem \eqref{DR_Dproblem} with boundary conditions \eqref{test1_DBCs}, initial condition \eqref{test1_ICs} and $16641$ spatial grid points ($129$ in each direction). The errors in the infinity norm are computed at $t=0.1$  by comparing the numerical solution to a reference solution, obtained with VODE applied to the unsplit problem.}
		\def\arraystretch{1.3}
		\begin{tabular}{C{1.5cm}C{1.5cm}C{1.5cm}C{0.4cm}C{1.5cm}C{1.5cm}}
			& \multicolumn{2}{c}{Strang} & & \multicolumn{2}{c}{modified Strang (exact elliptic)} \\
			\hline
			step size & $l^{\infty}$ error & order & & $l^{\infty}$ error & order \\
			\hline
			\ssp$2.50 \cdot 10^{-2}$    & \ssp$1.51\cdot 10^{-2}$    & \ssp --     & & \ssp$3.35\cdot 10^{-4}$    & \ssp --   \\
			$1.25 \cdot 10^{-2}$    & $7.45\cdot 10^{-3}$    & $1.02$ & & $8.28\cdot 10^{-5}$    & $2.02$  \\
			$6.25\cdot 10^{-3}$    & $3.66\cdot 10^{-3}$    & $1.03$ & & $2.08\cdot 10^{-5}$    & $1.99$  \\
			$3.13\cdot 10^{-3}$    & $1.78\cdot 10^{-3}$    & $1.04$ & & $5.24\cdot 10^{-6}$    & $1.99$  \\
			$1.56\cdot 10^{-3}$    & $8.47\cdot 10^{-4}$    & $1.07$ & & $1.36\cdot 10^{-6}$    & $1.94$  \\
		\end{tabular}
		\newline
		\newline
		\newline
		\begin{tabular}{C{1.5cm}C{1.5cm}C{1.5cm}C{0.003cm}C{1.5cm}C{1.5cm}}
			& \multicolumn{2}{c}{modified Strang ($ q_n = f(u_n) $)} & & \multicolumn{2}{c}{modified Strang (grid average)} \\
			\hline
			step size & $l^{\infty}$ error & order & &$l^{\infty}$ error & order \\
			\hline
			\ssp$2.50\cdot 10^{-2}$    & \ssp$6.49\cdot 10^{-4}$    & \ssp --    	& & \ssp$3.35\cdot 10^{-4}$    & \ssp --   \\
			$1.25\cdot 10^{-2}$    & $1.62\cdot 10^{-4}$    & $2.00$   & & $8.28\cdot 10^{-5}$    & $2.02$  \\
			$6.25\cdot 10^{-3}$    & $4.02\cdot 10^{-5}$    & $2.01$  & & $2.08\cdot 10^{-5}$    & $1.99$   \\
			$3.13\cdot 10^{-3}$    & $9.84\cdot 10^{-6}$    & $2.03$  & & $5.25\cdot 10^{-6}$    & $1.99$   \\
			$1.56\cdot 10^{-3}$    & $2.37\cdot 10^{-6}$    & $2.05$  & & $1.36\cdot 10^{-6}$    & $1.95$   \\
		\end{tabular}
	\label{Tab:D_InfNorm}
	\end{table}
\end{center}

\begin{center}
	\begin{table}[h!t!b!]
		\caption{Errors and observed orders of standard and modified Strang splitting schemes for the integration of problem \eqref{DR_Dproblem} with boundary conditions \eqref{test1_DBCs}, initial condition \eqref{test1_ICs} and $16641$ spatial grid points ($129$ in each direction). The errors in the $2$-norm are computed at $t=0.1$  by comparing the numerical solution to a reference solution, obtained with VODE applied to the unsplit problem.}
		\def\arraystretch{1.2}
		\begin{tabular}{C{1.5cm}C{1.5cm}C{1.5cm}C{0.4cm}C{1.5cm}C{1.5cm}}
			& \multicolumn{2}{c}{Strang} & & \multicolumn{2}{c}{modified Strang (exact elliptic)} \\
			\hline
			step size & $l^{2}$ error & order & & $l^{2}$ error & order \\
			\hline
			\ssp$2.50\cdot 10^{-2}$    &   \ssp$ 5.65 \cdot 10^{-3} $      & \ssp --     & &  \ssp$2.39\cdot 10^{-4}$    & \ssp --  \\
			$1.25\cdot 10^{-2}$    &   $ 2.31 \cdot 10^{-3} $      & $1.29$ & &  $6.06\cdot 10^{-5}$    & $1.98$ \\
			$6.25\cdot 10^{-3}$    &   $ 9.43 \cdot 10^{-4} $      & $1.30$ & &  $1.52\cdot 10^{-5}$    & $2.00$  \\
			$3.13\cdot 10^{-3}$    &   $ 3.80 \cdot 10^{-4} $      & $1.31$ & &  $3.81\cdot 10^{-6}$    & $2.00$  \\
			$1.56\cdot 10^{-3}$    &   $ 1.50 \cdot 10^{-4} $      & $1.34$ & &  $9.74\cdot 10^{-7}$    & $1.97$ \\
			
		\end{tabular}
		\newline
		\newline
		\newline
		\begin{tabular}{C{1.5cm}C{1.5cm}C{1.5cm}C{0.003cm}C{1.5cm}C{1.5cm}}
			& \multicolumn{2}{c}{modified Strang ($ q_n = f(u_n) $)} & & \multicolumn{2}{c}{modified Strang (grid average)} \\
			\hline
			step size & $l^{2}$ error & order & &$l^{2}$ error & order \\
			\hline
			\ssp$2.50\cdot 10^{-2}$    &   \ssp$ 2.33 \cdot 10^{-4}$    & \ssp --      & &  \ssp$ 2.39\cdot 10^{-4}$    & \ssp --   \\
			$1.25\cdot 10^{-2}$    &   $ 4.81 \cdot 10^{-5}$    & $2.27$  & &  $ 6.06\cdot 10^{-5}$    & $1.98$  \\
			$6.25\cdot 10^{-3}$    &   $ 1.05 \cdot 10^{-5}$    & $2.20$  & &  $ 1.52\cdot 10^{-5}$    & $2.00$   \\
			$3.13\cdot 10^{-3}$    &   $ 2.35 \cdot 10^{-6}$    & $2.16$  & &  $ 3.81\cdot 10^{-6}$    & $2.00$   \\
			$1.56\cdot 10^{-3}$    &   $ 5.39 \cdot 10^{-7}$    & $2.12$  & &  $ 9.74\cdot 10^{-7}$    & $1.97$  \\
			
		\end{tabular}
			\label{Tab:D_2Norm}
	\end{table}
\end{center}

\begin{figure}[h!t!b!]
	\centering
	\begin{tabular}{ll}
		\subfigure[Exact elliptic correction]
		{\includegraphics[width=5cm,height=5cm]{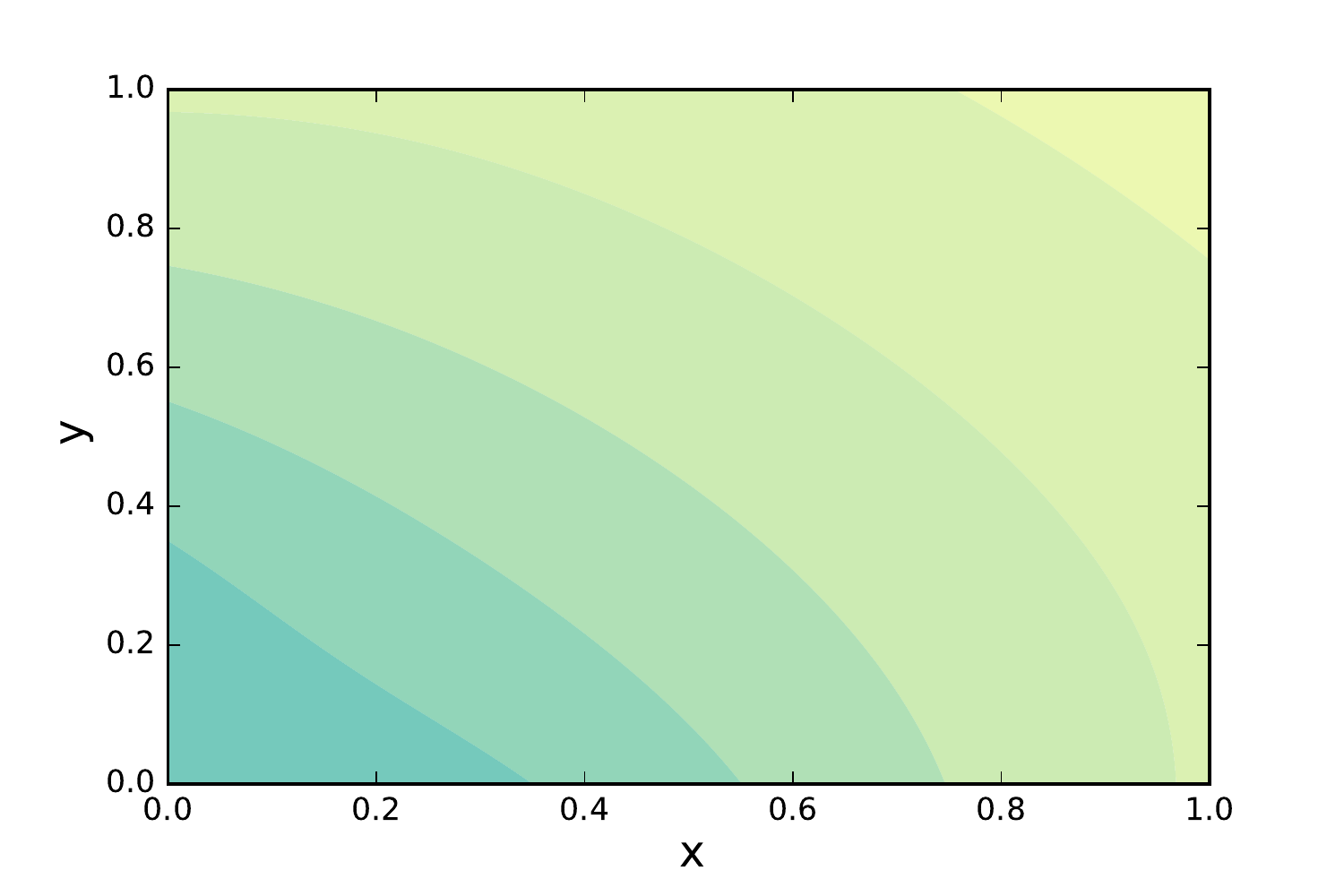}\label{Fig:DCorrFuncSolve}}&
		\subfigure[$q_n=f(u_n)$]
		{\includegraphics[width=5cm,height=5cm]{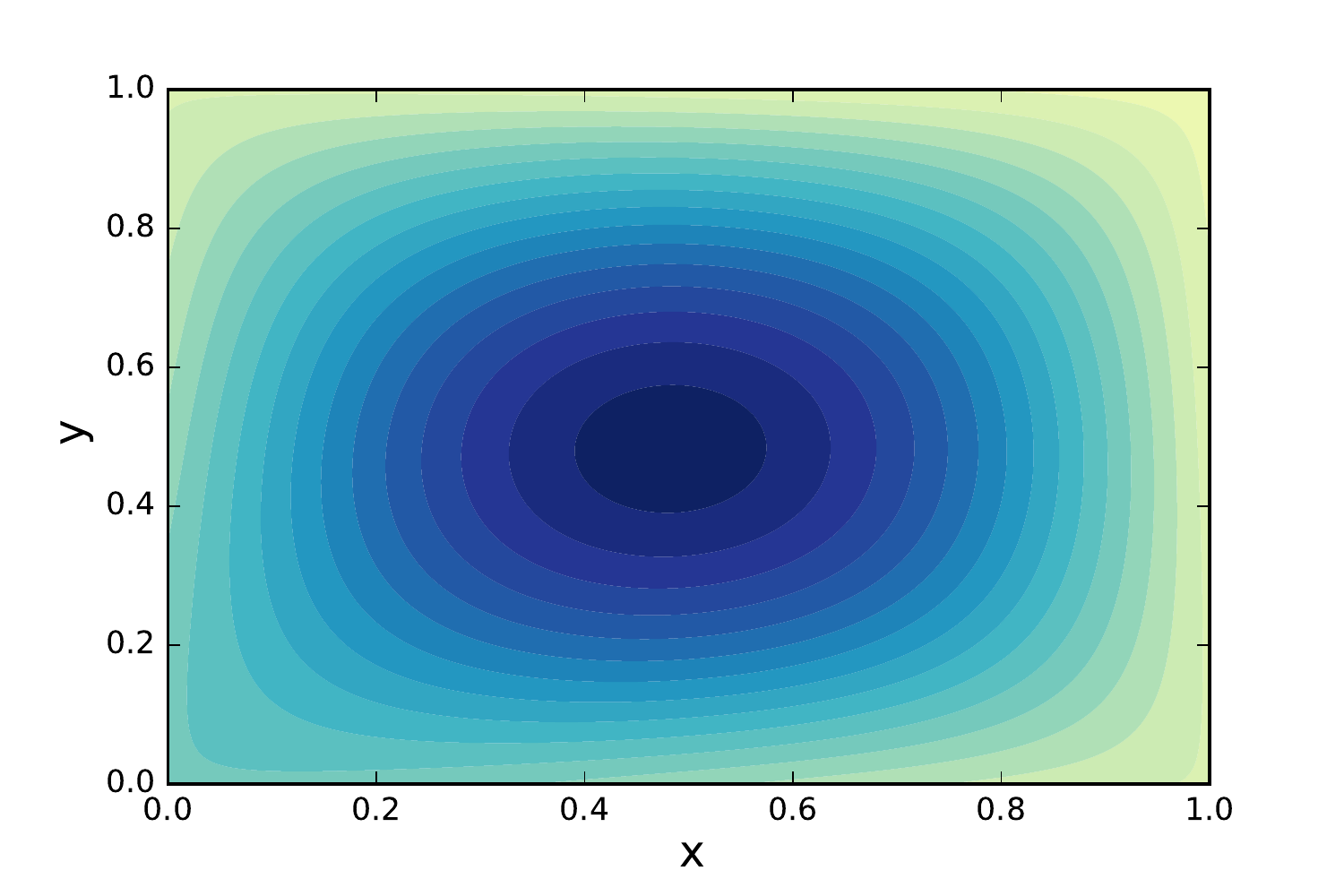}\label{Fig:DCorrFuncDirect}}\\
		\multicolumn{2}{c}{
		\subfigure[Grid average correction]
		{\includegraphics[width=6cm,height=5cm]{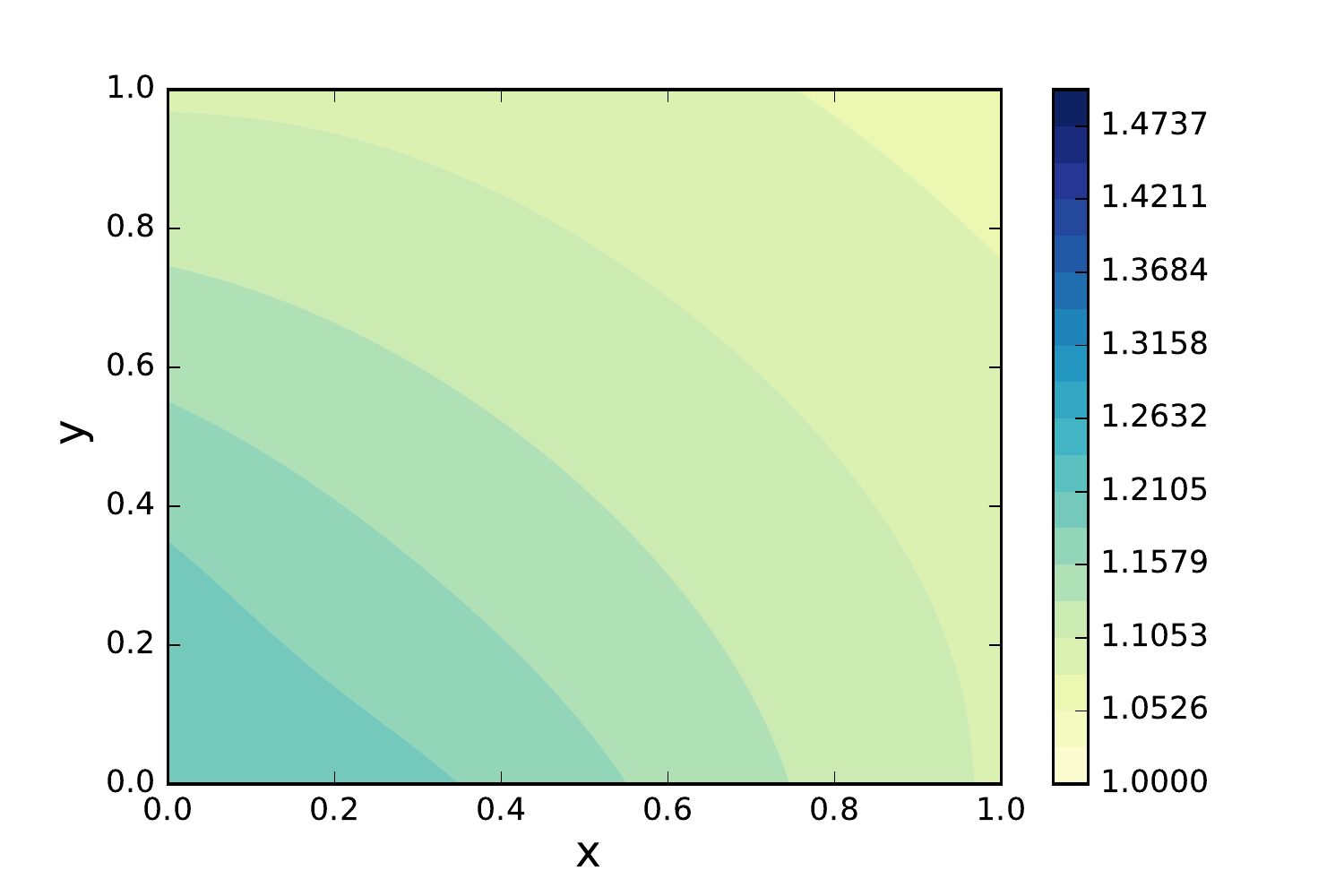}\label{Fig:DCorrFuncGrid}} }
	\end{tabular}
	\caption{Correction functions at $t=0.1$ constructed as the solution of the elliptic problem \eqref{SolveProblem} (top left), as $q_n=f(u_n)$ (top right) and as the average of neighbouring points (bottom) for the integration of problem \eqref{DR_Dproblem} with boundary conditions \eqref{test1_DBCs}, initial condition \eqref{test1_ICs} and 16641 spatial grid points (129 in each direction).} \label{Fig:DCorrFunc}
\end{figure}

\end{example}


\section{The Neumann case} \label{sec:Neumann}
In this section, we integrate system \eqref{system} subject to inhomogeneous Neumann boundary conditions
\begin{equation}\label{NBCs}
\partial_n u |_{\partial \Omega} = b.
\end{equation}
Following the idea presented in \cite{EOPart2}, we aim to construct a correction function $q$ such that
\begin{equation}
\partial_n q |_{\partial \Omega} = \partial_n f(u)|_{\partial \Omega}
\end{equation}
which is equivalent to
\begin{equation} \label{NCorrFunc}
\partial_n q |_{\partial \Omega} =f'(u)|_{\partial \Omega} \, b.
\end{equation}
In the following, we describe three different procedures to achieve this goal.

First, we can construct the correction function by solving the elliptic problem
\begin{subequations}\label{NSolveProblem}
\begin{equation}\label{NSolveProblem-a}
\begin{split}
\Delta q &= g, \\
\partial_n q |_{\partial \Omega} &=f'(u)|_{\partial \Omega} \, b,
\end{split}
\end{equation}
through, for instance, a direct method or a Krylov subspace scheme. Note that the constant $g$ is determined by the compatibility condition
\begin{equation}\label{NSolveProblem-constraint}
	g = \frac1{|\Omega|}\int_{\partial \Omega} f'(u) \, b \, \text{d} s.
\end{equation}
\end{subequations}

Second, we observe again that the choice $q_n=f(u_n)$ is compatible with the boundary condition \eqref{NCorrFunc} for the nonlinearity $f$. Thus, we can obtain the correction function as in the Dirichlet case.

Note that the moving average approach described in Algorithm \ref{algorithm:moving_average} for the Dirichlet case could change the slope of $u$ at the boundary. Therefore, we introduce as a third possibility an alternative technique based on a multigrid-like approach. In Section \ref{sec:Dirichlet}, we have already mentioned that multigrid schemes are widely employed to solve elliptic (and consequently parabolic) problems because they significantly reduce the number of iterations used. Since multigrid methods consist in relaxing the linear system stemming from the discretization of the considered PDE on a hierarchy of grids, they need some strategies to transfer information from a coarse grid to a fine grid (\emph{prolongation}) and, conversely, from a fine grid to a coarse grid (\emph{restriction}). The prolongation is commonly carried out by means of interpolation techniques. In this paper, we use linear interpolation. Let $V^H$ be a grid function, defined on the regular grid
\begin{equation}
\Omega^H = \{ (iH,jH)\mid 0\le i,j\le 2N\}
\end{equation}
with mesh width $H$. Interpolation from the coarse grid with mesh width $2H$ to the fine grid with width $H$ is achieved by
\begin{subequations} \label{2DInterp}
	\begin{align}
		V_{2i,2j}^H&=V_{i,j}^{2H}, \quad 0\le i,j\le N, \\
		V_{2i+1,2j}^H&=\dfrac{1}{2}\Big(V_{i,j}^{2H}+V_{i+1,j}^{2H}\Big), \quad 0\le i\le N-1, \ 0\le j\le N,\\
		V_{2i,2j+1}^H&=\dfrac{1}{2}\Big(V_{i,j}^{2H}+V_{i,j+1}^{2H}\Big), \quad 0\le i\le N, \ 0\le j\le N-1,
\shortintertext{and}
		V_{2i+1,2j+1}^H&=\dfrac{1}{4}\Big(V_{i,j}^{2H}+V_{i+1,j}^{2H}+V_{i,j+1}^{2H}+V_{i+1,j+1}^{2H}\Big)
	\end{align}
\end{subequations}
for $i,j=0, \dots, N-1$. In this procedure, the values of the two-dimensional array $V^H$ at even-numbered fine-grid points are equal to the values of $V^{2H}$  at coarse-grid points, whereas the values of $V^H$ at mixed and odd-numbered fine-grid points are computed as the average of the values of $V^{2H}$  at neighbouring coarse-grid points. On the other hand, the simplest restriction operator is the injection, which consists in mapping function values on the even-numbered fine-grid points to the coarse-grid points as follows
\begin{equation}\label{2Drestr}
	V_{i,j}^{2H} = V_{2i,2j}^H,   \quad i,j=0, \dots, N.
\end{equation}
In this work, we simply employ one half of a \textit{V-cycle}, as described in Algorithm~\ref{algorithm:multigrid} (which means that we do not need to make use of the restriction operator). For this purpose, we discretize the system \eqref{NSolveProblem} as $Aq=b$. As the basic smoothing step, we use weighted Jacobi iterations
\begin{equation}\label{Jacobi}
	q^{(s+1)}=(I- \omega D^{-1}A)q^{(s)} + \omega D^{-1}b,
\end{equation}
where $D$ is the diagonal of $A$ and $\omega$ is usually chosen equal to $2/3$. Another possibility are Gauss--Seidel iterations
\begin{equation}\label{GaussSeidel}
q^{(s+1)}=(D-L)^{-1}Uq^{(s)} + (D-L)^{-1}b,
\end{equation}
where $-L$ and $-U$ are the strictly lower and upper triangular parts of $A$, respectively.

\begin{algorithm}[t]
	\caption{How to compute the values of the correction function $q$ by employing one half of a V-cycle. }
	\label{algorithm:multigrid}
	\begin{enumerate}
        \vspace{1mm}
		\item Choose an integer $s \ge 1$, set $\ell=2^s$ and $H=\ell h/2$.
        \vspace{1mm}
		\item Solve the discretized system $Aq=b$ on the coarsest grid $\Omega^{2H}$ through a direct method to obtain $q^{2H}$.
		\vspace{1mm}
        \item Compute the prolongation $q^H$ by means of \eqref{2DInterp} and set $\ell=\ell/2$. \label{alg2:entry}
        \vspace{1mm}
		\item Relax the system $A^H q^H=b^H$ by using $\nu$ Jacobi or Gauss--Seidel iterations \eqref{Jacobi} with starting value $q^H$.
		\vspace{1mm}	
		\item If $\ell=1$ set $q=q^H$ and stop. Else set $H=H/2$ and go to \ref{alg2:entry}.
	\end{enumerate}
\end{algorithm}


\subsection{Numerical results for Neumann boundary conditions} \label{subsec:NeumannNumExp}
In this section, we discuss some numerical experiments for the following diffusion-reaction problem
\begin{subequations}\label{DR_Nproblem}
	\begin{alignat}{3}
	\partial_t u &= \Delta u + u^2,  \label{DR_Nequation} \\
	\partial_n u|_{\partial \Omega} &= b, \label{DR_NBCs}\\
	u(0)&=u_0 \label{DR_NICs}	
	\end{alignat}
\end{subequations}
on $\Omega=(0,1)^2$.
In the experiments, we employ the same methods mentioned in Section \ref{subsec:DirichletNumExp} to discretize the Laplacian and integrate each sub-problem of the splitted system \eqref{ModSplitSystem}.
In the remainder of this section, the terms \emph{Strang} and \emph{modified Strang} will refer to the standard Strang splitting and the modified Strang splitting schemes described in Section \ref{sec:problem}, respectively. Among the strategies of constructing the correction function proposed in Section \ref{sec:Neumann}, we refer to \emph{exact elliptic} as the solution of the Laplace problem endowed with Neumann boundary conditions \eqref{NSolveProblem} and to \emph{Jacobi} as the procedure presented in Algorithm \ref{algorithm:multigrid}. The computation of the correction function through the values of the numerical solution at the beginning of the current time step is again labelled as $q_n=f(u_n)$.

\begin{example}
We first integrate problem \eqref{DR_Nproblem} with the following boundary and initial data:
\begin{subequations}\label{N1}
	\begin{alignat}{2}
        b(x,y)&=\begin{cases}
        \phantom{-}1 & \ (x=1) \vee (y=1),\\
        -1 & \ (x=0)\vee (y=0),
        \end{cases}  \label{N1_BCs} \\
		u(0,x,y) &= x+y+ \sin^2(\pi \, x) \, \sin^2(\pi \, y),\qquad &(x,y)\in\Omega. \label{N1_ICs}
	\end{alignat}
\end{subequations}
Tables \ref{Tab:Nb1_InfNorm} and \ref{Tab:Nb1_2Norm} show that the standard Strang splitting exhibits orders $1.5$ and $1.75$ for the infinity and the $2$-norm, respectively, as is proved in \cite{EOPart2}. On the other hand, the modification presented in Section \ref{sec:problem} combined with the strategies of constructing the correction function proposed in Section \ref{sec:Neumann} leads to Strang splitting schemes of order $2$. Therefore, the problem of order reduction arises even with these simple Neumann boundary conditions for the standard Strang splitting, but the modifications described in Section \ref{sec:Neumann} allow us to overcome this limitation. Moreover, all the modified Strang splitting schemes exhibit higher accuracy than the standard Strang splitting. These favourable results can also be observed in more challenging problems, as reported in the following examples.

\begin{center}
	\begin{table}[t!h!b!]
		\caption{Errors and observed orders of standard and modified Strang splitting schemes for the integration of problem \eqref{DR_Nproblem} with boundary conditions \eqref{N1_BCs}, initial condition \eqref{N1_ICs} and $16641$ spatial grid points ($129$ in each direction). The errors in the infinity norm are computed at $t=0.1$  by comparing the numerical solution to a reference solution, obtained with VODE applied to the unsplit problem. For Jacobi Strang, three Jacobi iterations have been used.}
		\def\arraystretch{1.2}
		\begin{tabular}{C{1.5cm}C{1.5cm}C{1.5cm}C{0.5cm}C{1.5cm}C{1.5cm}}
			 & \multicolumn{2}{c}{Strang} & & \multicolumn{2}{c}{modified Strang (exact elliptic)} \\
			\hline
			step size & $l^{\infty}$ error & order & & $l^{\infty}$ error & order \\
			\hline
            \msp$2.50 \formExp{-2}$    & \msp$9.40\formExp{-3}$    & \msp --      & & \msp$2.37\formExp{-3}$    & \msp --    \\
			$1.25 \formExp{-2}$    & $3.38\formExp{-3}$    & $1.48$  & & $5.62\formExp{-4}$    & $2.08$  \\
			$6.25 \formExp{-3}$    & $1.20\formExp{-3}$    & $1.49$  & & $1.34\formExp{-4}$    & $2.07$    \\
			$3.13 \formExp{-3}$    & $4.26\formExp{-4}$    & $1.50$  & & $3.24\formExp{-5}$    & $2.05$      \\
			$1.56 \formExp{-3}$    & $1.50\formExp{-4}$    & $1.50$  & & $7.88\formExp{-6}$    & $2.04$   \\
		\end{tabular}
		\newline
		\newline
		\newline
		\begin{tabular}{C{1.5cm}C{1.5cm}C{1.5cm}C{0.05cm}C{1.5cm}C{1.5cm}}
			& \multicolumn{2}{c}{modified Strang ($q_n=f(u_n)$)} & & \multicolumn{2}{c}{modified Strang (Jacobi)} \\
			\hline
			step size & $l^{\infty}$ error & order & &$l^{\infty}$ error & order \\
			\hline
			\msp$2.50\formExp{-2}$    & \msp$3.13\formExp{-3}$    & \msp --    	& & \msp$2.43\formExp{-3}$    & \msp --     \\
			$1.25\formExp{-2}$    & $7.52\formExp{-4}$    & $2.06$  & & $5.54\formExp{-4}$    & $2.08$  \\
			$6.25\formExp{-3}$    & $1.82\formExp{-4}$    & $2.05$  & & $1.32\formExp{-4}$    & $2.07$  \\
			$3.13\formExp{-3}$    & $4.44\formExp{-5}$    & $2.04$  & & $3.19\formExp{-5}$    & $2.05$   \\
			$1.56\formExp{-3}$    & $1.09\formExp{-5}$    & $2.03$  & & $7.79\formExp{-6}$    & $2.03$ \\
		\end{tabular}
		\label{Tab:Nb1_InfNorm}
	\end{table}
\end{center}

\begin{center}
	\begin{table}[t!h!b!]
		\caption{Errors and observed orders of standard and modified Strang splitting schemes for the integration of problem \eqref{DR_Nproblem} with boundary conditions \eqref{N1_BCs}, initial condition \eqref{N1_ICs} and $16641$ spatial grid points ($129$ in each direction). The errors in the $2$-norm are computed at $t=0.1$  by comparing the numerical solution to a reference solution, obtained with VODE applied to the unsplit problem. For Jacobi Strang, three Jacobi iterations have been used.}
		\def\arraystretch{1.2}
		\begin{tabular}{C{1.5cm}C{1.5cm}C{1.5cm}C{0.5cm}C{1.5cm}C{1.5cm}}
			& \multicolumn{2}{c}{Strang} & & \multicolumn{2}{c}{modified Strang (exact elliptic)} \\
			\hline
			step size & $l^{2}$ error & order & & $l^{2}$ error & order \\
			\hline
			\msp$2.50\formExp{-2}$    & \msp$1.31\formExp{-5}$    & \msp --      & & \msp$ 7.72 \formExp{-6}$ & \msp -- \\
			$1.25\formExp{-2}$    & $3.92\formExp{-6}$    & $1.75$  & & $ 2.03 \formExp{-7}$ & $1.93 $  \\
			$6.25\formExp{-3}$    & $1.17\formExp{-6}$    & $1.74$  & & $ 5.16 \formExp{-7}$ & $ 1.97 $  \\
			$3.13\formExp{-3}$    & $3.54\formExp{-7}$    & $1.72$  & & $1.30 \formExp{-7}$ & $1.99 $  \\
			$1.56\formExp{-3}$    & $1.09\formExp{-7}$    & $1.70$  & & $ 3.31 \formExp{-8}$ & $1.97$  \\
		
		\end{tabular}
		\newline
		\newline
		\newline
		\begin{tabular}{C{1.5cm}C{1.5cm}C{1.5cm}C{0.05cm}C{1.5cm}C{1.5cm}}
			& \multicolumn{2}{c}{modified Strang ($q_n=f(u_n)$)} & & \multicolumn{2}{c}{modified Strang (Jacobi)} \\
			\hline
			step size & $l^{2}$ error & order & &$l^{2}$ error & order \\
			\hline
			\msp$2.50\formExp{-2}$    & \msp$9.26\formExp{-6}$    & \msp--      & & \msp$7.67 \formExp{-6}$ & \msp -- \\
			$1.25\formExp{-2}$    & $2.38\formExp{-6}$    & $1.96$  & & $ 2.01\formExp{-6}$  &  $1.93$ \\
			$6.25\formExp{-3}$    & $6.03\formExp{-7}$    & $1.98$  & & $ 5.12\formExp{-7}$ & $1.97$\\
			$3.13\formExp{-3}$    & $1.52\formExp{-7}$    & $1.99$  & & $1.29\formExp{-7}$  & $ 1.99$\\
			$1.56\formExp{-3}$    & $3.82\formExp{-8}$    & $1.99$  & & $3.29 \formExp{-8}$  & $ 1.97$			
		\end{tabular}
		\label{Tab:Nb1_2Norm}
	\end{table}
\end{center}
\end{example}

\begin{example}
Next, we examine problem \eqref{DR_Nequation} with initial condition
\begin{equation}\label{N3_ICs}
u(0,x,y) = 3+ \mathrm{e}^{-5(y-0.5)^2} \cos(2\pi(x+y))
\end{equation}
and Neumann boundary data
\begin{equation}\label{N5_BCs}
b(t,x,y)=
\begin{cases}
    (-1)^x 2\pi \,\mathrm{e}^{-5(y-0.5)^2} \sin \theta, & x\in \set{0,1}, \;  y \in [0,1], \\[3mm]
    (-1)^y 2\,\mathrm{e}^{-5(y-0.5)^2}\left(5(y-0.5)\cos \theta+\pi \sin \theta \right), &y\in \set{0,1}, \; x \in [0,1],
\end{cases}
\end{equation}
where $\theta=2 \pi(x+y)$. As shown in Tables \ref{Tab:Nb2_InfNorm} and \ref{Tab:Nb2_2Norm}, the standard Strang splitting has orders $1.5$ and $1.75$ for the infinity and the $2$-norm, respectively, as it is proved in \cite{EOPart2}. The modified schemes described in Section \ref{sec:Neumann}, however, have order two and higher accuracy.
\end{example}

\begin{center} 
	\begin{table}[t]
		\caption{Errors and orders of standard and modified Strang splitting schemes for the integration of problem \eqref{DR_Nproblem} with initial condition \eqref{N3_ICs}, boundary conditions \eqref{N5_BCs} and $16641$ spatial grid points ($129$ in each direction). The errors in the infinity norm are computed at $t=0.1$  by comparing the numerical solution to a reference solution, obtained with VODE applied to the unsplit problem. For Jacobi Strang, three Jacobi iterations have been used.}
		\def\arraystretch{1.3}
		\begin{tabular}{C{1.5cm}C{1.5cm}C{1.5cm}C{0.5cm}C{1.5cm}C{1.5cm}}
			& \multicolumn{2}{c}{Strang} & & \multicolumn{2}{c}{modified Strang (exact elliptic)} \\
			\hline
			step size & $l^{\infty}$ error & order & & $l^{\infty}$ error & order \\
			\hline
			\msp$2.50 \formExp{-2}$    & \msp$4.84\formExp{-2}$    & \msp --      & & \msp$1.24 \formExp{-2}$    & \msp --    \\
			$1.25 \formExp{-2}$    & $1.75 \formExp{-2}$    & $1.47$  & & $3.18 \formExp{-3}$    & $1.96$  \\
			$6.25 \formExp{-3}$    & $6.19 \formExp{-3}$    & $1.50$  & & $7.91 \formExp{-4}$    & $2.01$   \\
			$3.13 \formExp{-3}$    & $2.18 \formExp{-3}$    & $1.51$  & & $1.95\formExp{-4}$    & $2.02$  \\
			$1.56 \formExp{-3}$    & $7.61 \formExp{-4}$    & $1.52$  & & $4.81 \formExp{-5}$    & $2.02$  \\
		\end{tabular}
		\newline
		\newline
		\newline
		\begin{tabular}{C{1.5cm}C{1.5cm}C{1.5cm}C{0.05cm}C{1.5cm}C{1.5cm}}
			& \multicolumn{2}{c}{modified Strang ($q_n=f(u_n)$)} & & \multicolumn{2}{c}{modified Strang (Jacobi)} \\
			\hline
			step size & $l^{\infty}$ error & order & &$l^{\infty}$ error & order \\
			\hline
			\msp$2.50 \formExp{-2}$    & \msp$1.45 \formExp{-2}$    & \msp --       & & \msp$1.33 \formExp{-2}$    & \msp --    \\
			$1.25 \formExp{-2}$    & $3.74 \formExp{-3}$    & $1.96$   & & $3.50\formExp{-3}$    & $1.92$  \\
			$6.25 \formExp{-3}$    & $9.35 \formExp{-4}$    & $2.00$    & & $8.95 \formExp{-4}$    & $1.97$   \\
			$3.13 \formExp{-3}$    & $2.31 \formExp{-4}$    & $2.02$   & & $2.27 \formExp{-4}$    & $1.98$  \\
			$1.56 \formExp{-3}$    & $5.70 \formExp{-5}$    & $2.02$   & & $5.78 \formExp{-5}$    & $1.97$   \\
			
		\end{tabular}
		\label{Tab:Nb2_InfNorm}
	\end{table}
\end{center}

\subsection{Numerical results for mixed boundary conditions}
We now present a more challenging problem, involving mixed boundary conditions. In particular, we consider system \eqref{DR_Nequation} on $\Omega=(0,1)^2$, provided with initial conditions \eqref{DR_NICs}, Dirichlet boundary conditions on the top and bottom of the domain $\Omega$ and Neumann boundary conditions on the left and right side of $\Omega$. We thus consider the problem:
\begin{subequations}\label{mixBproblem}
	\begin{alignat}{3}
	\partial_t u &= \Delta u + u^2,  \label{mixBequation} \\
	\partial_n u|_{\partial \Omega_N} &= b_N, \quad \partial \Omega_N=\set{x\in \{0,1\}, \;y\in[0,1]},   \label{mixN}\\
	u|_{\partial \Omega_D} &= b_D, \quad \partial \Omega_D=\set{x\in [0,1] , \;y\in\{0,1\}},   \label{mixD}\\
	u(0)&=u_0 \label{mixICs}.	
	\end{alignat}
\end{subequations}
We note that the more complicated the boundary conditions are, the more Jacobi iterations are needed to guarantee the smoothness of the correction function. Therefore, it may be more suitable to use Gauss--Seidel instead. Henceforth, we refer to \emph{Gauss--Seidel} or \emph{Jacobi} as the method used in Algorithm~\ref{algorithm:multigrid}, according to the kind of smoother used.

\begin{center}
	\begin{table}[t]
		\caption{Errors and orders of standard and modified Strang splitting schemes for the integration of problem \eqref{DR_Nproblem} with initial condition \eqref{N3_ICs}, boundary conditions \eqref{N5_BCs} and $16641$ spatial grid points ($129$ in each direction). The errors in the $2$-norm are computed at $t=0.1$  by comparing the numerical solution to a reference solution, obtained with VODE applied to the unsplit problem. For Jacobi Strang, three Jacobi iterations have been used.}
		\def\arraystretch{1.3}
		\begin{tabular}{C{1.5cm}C{1.5cm}C{1.5cm}C{0.5cm}C{1.5cm}C{1.5cm}}
			& \multicolumn{2}{c}{Strang} & & \multicolumn{2}{c}{modified Strang (exact elliptic)} \\
			\hline
			step size & $l^{2}$ error & order & & $l^{2}$ error & order \\
			\hline
			\msp$2.50 \formExp{-2}$    & \msp$6.44\formExp{-5} $    & \msp --      & & \msp$3.45\formExp{-5}$     & \msp --  \\
			$1.25 \formExp{-2}$    & $2.03\formExp{-5} $    & $1.67$  & & $9.90\formExp{-6}$     & $1.80$  \\
			$6.25 \formExp{-3}$    & $6.18 \formExp{-6}$    & $1.71$  & & $6.74 \formExp{-7}$    & $1.96$  \\
			$3.13 \formExp{-3}$    & $1.88 \formExp{-6}$    & $1.71$  & & $6.74 \formExp{-7}$    & $1.96$   \\
			$1.56 \formExp{-3}$    & $5.81 \formExp{-7}$    & $1.70$  & & $ 1.71 \formExp{-7}$    & $1.98$  \\
			
		\end{tabular}
		\newline
		\newline
		\newline
		\begin{tabular}{C{1.5cm}C{1.5cm}C{1.5cm}C{0.05cm}C{1.5cm}C{1.5cm}}
			& \multicolumn{2}{c}{modified Strang ($q_n=f(u_n)$)} & & \multicolumn{2}{c}{modified Strang (Jacobi)} \\
			\hline
			step size & $l^{2}$ error & order & &$l^{2}$ error & order \\
			\hline
			\msp$2.50 \formExp{-2}$    & \msp$3.39 \formExp{-5}$    & \msp --     & & \msp $3.49 \formExp{-5}$    & \msp --      \\
			$1.25 \formExp{-2}$    & $9.69 \formExp{-6}$    & $1.90$  & & $1.00 \formExp{-5}$   & $1.80$    \\
			$6.25 \formExp{-3}$    & $2.59 \formExp{-6}$    & $1.95$ & & $2.66 \formExp{-6}$    & $1.91$    \\
			$3.13 \formExp{-3}$    & $6.67 \formExp{-7}$    & $1.98$ & & $6.82 \formExp{-7}$    & $1.96$   \\
			$1.56 \formExp{-3}$    & $1.69 \formExp{-7}$    & $1.99$ & & $1.73 \formExp{-7}$    & $1.98$
		\end{tabular}
		\label{Tab:Nb2_2Norm}
	\end{table}
\end{center}

\begin{example}
We integrate problem \eqref{mixBequation} provided with the initial condition
\begin{equation}\label{mix_ICs}
u(0,x,y) = 3+ \mathrm{e}^{-5(y-0.5)^2} \cos(2\pi(x+y))
\end{equation}
and the following boundary data
\begin{equation}\label{mix_BCs}
\begin{split}
b_N(t,x,y) & = (-1)^x 2\pi \,\mathrm{e}^{-5(y-0.5)^2} \sin(2\pi(x+y)), \\
b_D(t,x,y) & = 3+ \mathrm{e}^{-5/4} \cos(2\pi(x+y)).
\end{split}
\end{equation}
Tables \ref{Tab:mix_InfNorm} and \ref{Tab:mix_2Norm} show that the standard Strang splitting has orders $1$ and $1.3$ for the infinity and the $2$-norm, respectively. This result is consistent with the results obtained in \cite{EOPart1}, since, in case of mixed boundary conditions, the order reduction due to the Dirichlet boundary conditions is dominant. The modified schemes described in Section \ref{sec:Neumann}, however, exhibit order two and much higher accuracy. Moreover, Table~\ref{Tab:mix_InfNorm_JvsGS} shows that the modified Strang splitting with Gauss--Seidel iterations is slightly more accurate and has an order closer to 2 than the corresponding method with Jacobi iterations.

\begin{center}
\begin{table}[t!h!b!]
	\caption{Errors and orders of standard and modified Strang splitting schemes for the integration of problem \eqref{mixBequation} with initial condition \eqref{mix_ICs}, mixed boundary conditions \eqref{mix_BCs} and $16641$ spatial grid points ($129$ in each direction). The errors in the infinity norm are computed at $t=0.1$ by comparing the numerical solution to a reference solution, obtained with VODE applied to the unsplit problem. For Gauss--Seidel Strang, twenty Gauss--Seidel iterations have been used.}
	\def\arraystretch{1.3}
	\begin{tabular}{C{1.5cm}C{1.5cm}C{1.5cm}C{0.5cm}C{1.5cm}C{1.5cm}}
		& \multicolumn{2}{c}{Strang} & & \multicolumn{2}{c}{modified Strang (exact elliptic)} \\
		\hline
		step size & $l^{\infty}$ error & order & & $l^{\infty}$ error & order \\
		\hline
		\msp$2.50 \formExp{-2}$    & \msp$1.31 \formExp{-1}$    & \msp --      & & \msp$ 1.16\formExp{-2}$    & \msp--    \\
		$1.25 \formExp{-2}$    & $ 6.15\formExp{-2}$    & $1.09$  & & $ 2.91\formExp{-3}$    & $2.00$  \\
		$6.25 \formExp{-3}$    & $ 2.89\formExp{-2}$    & $1.09$  & & $ 7.26\formExp{-4}$    & $2.00$   \\
		$3.13 \formExp{-3}$    & $ 1.34\formExp{-2}$    & $1.11$  & & $ 1.83\formExp{-4}$    & $1.99$   \\
		$1.56 \formExp{-3}$    & $ 6.02\formExp{-3}$    & $1.15$  & & $4.67 \formExp{-5}$    & $1.97$  \\
	\end{tabular}
	\newline
	\newline
	\newline
	\begin{tabular}{C{1.5cm}C{1.5cm}C{1.5cm}C{0.05cm}C{1.5cm}C{1.5cm}}
		& \multicolumn{2}{c}{modified Strang ($q_n=f(u_n)$)} & & \multicolumn{2}{c}{modified Strang (Gauss--Seidel)} \\
		\hline
		step size & $l^{\infty}$ error & order & &$l^{\infty}$ error & order \\
		\hline
		\msp$2.50 \formExp{-2}$    & \msp$1.30\formExp{-2}$     &\msp --      & & \msp$ 1.53 \formExp{-2}$    & \msp --       \\
		$1.25 \formExp{-2}$    & $ 3.03\formExp{-3}$    & $2.11$  & & $ 4.00 \formExp{-3}$    & $ 1.93$      \\
		$6.25 \formExp{-3}$    & $7.17\formExp{-4}$     & $2.08$  & & $ 1.03 \formExp{-3}$    & $ 1.95$     \\
		$3.13 \formExp{-3}$    & $1.72\formExp{-4}$     & $2.06$  & & $ 2.68 \formExp{-4}$    & $1.94$     \\
		$1.56 \formExp{-3}$    & $4.19\formExp{-5}$    & $2.04$  & & $ 6.96 \formExp{-5}$    & $1.95$     \\
		
	\end{tabular}
	\label{Tab:mix_InfNorm}
\end{table}
\end{center}

\begin{center}
\begin{table}[t!h!b!]
\caption{
Errors and orders of standard and modified Strang splitting schemes for the integration of problem \eqref{mixBequation} with initial condition \eqref{mix_ICs}, mixed boundary conditions \eqref{mix_BCs} and $16641$ spatial grid points ($129$ in each direction). The errors in the $2$-norm are computed at $t=0.1$  by comparing the numerical solution to a reference solution, obtained with VODE applied to the unsplit problem. For Gauss--Seidel Strang, three Gauss--Seidel iterations have been used.}
\def\arraystretch{1.3}

\begin{tabular}{C{1.5cm}C{1.5cm}C{1.5cm}C{0.5cm}C{1.5cm}C{1.5cm}}
	& \multicolumn{2}{c}{Strang} & & \multicolumn{2}{c}{modified Strang (exact elliptic)} \\
	\hline
	step size & $l^{2}$ error & order & & $l^{2}$ error & order \\
	\hline
	\msp$2.50 \formExp{-2}$    & \msp$ 2.69 \formExp{-4}$    & \msp --   & & \msp$ 3.39 \formExp{-5}$    & \msp --    \\
	$1.25 \formExp{-2}$    & $ 1.06\formExp{-4}$    & $1.34$ & & $ 9.06\formExp{-6}$    & $1.90$     \\
	$6.25 \formExp{-3}$    & $ 4.23\formExp{-5}$    & $1.32$ & & $ 2.34\formExp{-6}$    & $1.95$    \\
	$3.13 \formExp{-3}$    & $ 1.69\formExp{-5}$    & $1.33$ & & $ 5.94\formExp{-7}$    & $1.98$     \\
	$1.56 \formExp{-3}$    & $ 6.64\formExp{-6}$    & $1.35$ & & $ 1.50\formExp{-7}$    & $1.99$    \\
\end{tabular}
\newline
\newline
\newline
\begin{tabular}{C{1.5cm}C{1.5cm}C{1.5cm}C{0.1cm}C{1.5cm}C{1.5cm}}
	& \multicolumn{2}{c}{modified Strang ($q_n=f(u_n)$)} & & \multicolumn{2}{c}{modified Strang (Gauss-Seidel)} \\
	\hline
	step size & $l^{2}$ error & order & &$l^{2}$ error & order \\
	\hline
	\msp$2.50 \formExp{-2}$    & \msp$ 3.00\formExp{-5}$    & \msp --      & & \msp$ 3.71 \formExp{-5}$    & \msp--  \\
	$1.25 \formExp{-2}$    & $ 7.76\formExp{-6}$    & $1.95$  & & $  9.92 \formExp{-6}$    & $1.91$  \\
	$6.25 \formExp{-3}$    & $ 1.99\formExp{-6}$    & $1.97$  & & $ 2.56 \formExp{-6}$    & $1.96$  \\
	$3.13 \formExp{-3}$    & $ 5.02\formExp{-7}$    & $1.98$  & & $ 6.49 \formExp{-7}$    & $1.98$  \\
	$1.56 \formExp{-3}$    & $1.26 \formExp{-7}$    & $1.99$   & & $ 1.64 \formExp{-7}$    & $1.99$  \\
	
\end{tabular}
\label{Tab:mix_2Norm}
\end{table}
\end{center}

\begin{center}
\begin{table}[t!h!b!]
\caption{Errors and orders of the modified Strang splitting schemes with Gauss--Seidel and Jacobi smoother (with $n$ iterations) for the integration of problem \eqref{mixBequation} with initial condition \eqref{mix_ICs}, mixed boundary conditions \eqref{mix_BCs} and $16641$ spatial grid points ($129$ in each direction). The errors in the infinity norm are computed at $t=0.1$  by comparing the numerical solution to a reference solution, obtained with VODE applied to the unsplit problem. }
\def\arraystretch{1.3}
\begin{tabular}{C{1.5cm}C{1.5cm}C{1.5cm}C{0.5cm}C{1.5cm}C{1.5cm}}
$n=3$ & \multicolumn{2}{c}{modified Strang (Gauss-Seidel)} & & \multicolumn{2}{c}{modified Strang (Jacobi)} \\
\hline
step size & $l^{\infty}$ error & order & & $l^{\infty}$ error & order \\
\hline
\msp$2.50 \formExp{-2}$    & \msp$ 2.34\formExp{-2}$    & \msp --      & & \msp$ 2.96 \formExp{-2}$    & \msp--  \\
$1.25 \formExp{-2}$    & $6.64\formExp{-3}$     & $1.82$  & & $ 8.99 \formExp{-3}$    & $1.72$  \\
$6.25 \formExp{-3}$    & $ 1.88\formExp{-3}$    & $1.82$  & & $2.65 \formExp{-3}$    & $1.76 $  \\
$3.13 \formExp{-3}$    & $ 5.21\formExp{-4}$    & $1.85$  & & $7.60 \formExp{-4}$    & $1.80$  \\
$1.56 \formExp{-3}$    & $1.42 \formExp{-4}$    & $1.87$   & & $ 2.12 \formExp{-4}$    & $1.84$  \\

\end{tabular}
\newline
\newline
\newline
\begin{tabular}{C{1.5cm}C{1.5cm}C{1.5cm}C{0.5cm}C{1.5cm}C{1.5cm}}
$n=5$ & \multicolumn{2}{c}{modified Strang (Gauss-Seidel)} & & \multicolumn{2}{c}{modified Strang (Jacobi)} \\
\hline
step size & $l^{\infty}$ error & order & & $l^{\infty}$ error & order \\
\hline
\msp$2.50 \formExp{-2}$    & \msp$2.05\formExp{-2}$    & \msp --      & & \msp$ 2.64 \formExp{-2}$    & \msp --  \\
$1.25 \formExp{-2}$    & $5.67\formExp{-3}$    & $1.86$  & & $ 7.78 \formExp{-3}$    & $1.76$  \\
$6.25 \formExp{-3}$    & $1.57\formExp{-3}$    & $1.85$  & & $2.22 \formExp{-3}$    & $1.81$  \\
$3.13 \formExp{-3}$    & $4.27\formExp{-4}$    & $1.88$  & & $6.17 \formExp{-4}$    & $1.85$  \\
$1.56 \formExp{-3}$    & $1.14 \formExp{-4}$    & $1.90$   & & $ 1.67 \formExp{-4}$    & $1.88$  \\
	
\end{tabular}
\newline
\newline
\newline
\begin{tabular}{C{1.5cm}C{1.5cm}C{1.5cm}C{0.5cm}C{1.5cm}C{1.5cm}}
$n=20$ & \multicolumn{2}{c}{modified Strang (Gauss-Seidel)} & & \multicolumn{2}{c}{modified Strang (Jacobi)} \\
\hline
step size & $l^{\infty}$ error & order & & $l^{\infty}$ error & order \\
\hline
\msp$2.50 \formExp{-2}$    & \msp$1.53\formExp{-2}$    & \msp --      & & \msp$ 1.91 \formExp{-2}$    & \msp --  \\
$1.25 \formExp{-2}$    & $ 4.00\formExp{-3}$    & $1.93 $  & & $ 5.26 \formExp{-3}$    & $1.86$  \\
$6.25 \formExp{-3}$    & $1.03\formExp{-3}$    & $1.95$  & & $1.41 \formExp{-3}$    & $1.90$  \\
$3.13 \formExp{-3}$    & $2.68\formExp{-4}$    & $1.94$  & & $3.70\formExp{-4}$    & $1.93 $  \\
$1.56 \formExp{-3}$    & $6.96 \formExp{-5}$    & $1.95$   & & $ 9.57 \formExp{-5}$    & $1.95$  \\
\end{tabular}
\label{Tab:mix_InfNorm_JvsGS}
\end{table}
\end{center}

\end{example}


\section{Conclusions} \label{sec:Conclusions}

Strang splitting for diffusion-reaction equations with nontrivial boundary conditions is prone to order reduction and requires an appropriate boundary correction. This can be achieved by a smooth function that satisfies the boundary conditions of the nonlinearity. In this paper, we discussed three different approaches to compute this correction. They all have different characteristics with respect to ease of implementation and efficiency.

The first approach relies on an additional solution of an elliptic problem, similar to the one that is solved in the integrator itself. Even for time-dependent boundary conditions, the computational overhead of this correction is moderate. The second one employs the computed numerical solution. Although this approach is conceptually very simple, it can be expensive if the evaluation of the nonlinear vector field is more expensive than the solution of the diffusion equation. The third proposed approach is based on iterative schemes and their smoothing properties. Its efficient implementation (with cost proportional to the number of unknowns) requires a hierarchy of grids in a way similar to multigrid methods.

All these corrections have their advantages and disadvantages depending on the properties of the evolution equation.

\bibliographystyle{plain}
\bibliography{bibliography}

\end{document}